# ANTICONTROL TECHNIQUES FOR SYSTEMS OF LORENZ TYPE

D. CONSTANTINESCU[1]

[1] University of Craiova, A. I. Cuza Street, 200585 Craiova, Romania
E-mail: dconsta@yahoo.com

The paper investigates some basic dynamical properties of a general system obtained from the Lorenz system using a non-linear feedback controller. We focus on the bifurcation of the equilibrium points and on the existence and the description of homoclinic and heteroclinic orbits. We present necessary conditions for the anticontrol of chaos in the considered system.

## 1. INTRODUCTION

Over the last decades chaos has gradually moved from simply being a curious phenomenon to one with practical significance and applications. The control of chaos (the stabilization of the chaotic orbits) or the anticontrol of chaos (the creation of chaos in a stable system) are key issues in applications where chaos is dangerous (engineering, comunications, many physical, chemical, biological, social phenomena) respectively important and useful (thorough liquid mixing with low power consumption, high-performance circuit design for tele-communication, collapse prevention of power systems, biomedical engineering applications to the human brain and heart).

In order to control or anticontrol the chaos a feedback controller can be used. We will exemplify this technique using a system obtained by applying a simple and implementable controller to the Lorenz system

$$x' = a(y - z), \quad y' = cx - xz - y, \quad z' = -bz + xy. \tag{1}$$

This system, introduced in [1] in 1963, is considered a pattern of the chaotic behaviour. For some values of the parameters (for example $a = 10$, $b = 8/3$, $c = 28$) it displays a very complex dynamics characterized by the sensitive dependence on initial conditions and the existence of a strange attractor [1]. Hundreds of interesting papers were devoted to the study of its dynamical properties and to the description of its attractor. For other values of the parameters (for example $a > 0$, $b > 0$, $c \in (0, 1)$) the system has a simple behaviour because the equilibrium point $O(0, 0, 0)$ is a global attractor. It is an ideal example to apply the control/anticontrol techniques.

The controlled (or anticontrolled) system considered in this paper is

$$\begin{cases} x' = a(y-x) \\ y' = cx - xz - y + u(x,y,z) = cx - xz - y + Mx + Ny + Pxz. \\ z' = -bz + xy \end{cases} \qquad (2)$$

It is obtained from (1) using the control term $u = Mx + Ny + Pxz$.

At least three particular cases of (2) are important both from theoretical and practical point of view and were extensively studied.

The Chen system

$$x' = a(y-x), \quad y' = (c-a)x + cy - xz, \quad z' = -bz + xy \qquad (3)$$

introduced in [2] has an attractor which is not topologically equivalent with the Loreanz' attractor. Some global bifurcation in Chen system were pointed out and anticontrol techniques were applied in [3], many interesting properties of the chaotic attractor were studied in [4], local bifurcation were studied in [5]. Chaos control was investigated in [6] existence of homoclinic and heteroclinic orbita was investigate in [7].

A transition between the Lorenz and the Chen system is the Lu system

$$x' = a(y-x), \quad y' = cy - xz, \quad z' = -bz + xy. \qquad (4)$$

introduced in [8]. The comparison of some dynamical properties of Lorenz, Chen and Lu systems is showed in [9].

The T-system (5) was introduced recently in [10].

$$x' = a(y-x), \quad y' = (c-a)x - axz, \quad z' = -bz + xy \qquad (5)$$

These systems have specific dynamical properties, but they have some common characteristics, due to the fact that they are particular cases of (2).

Indeed, the Lorenz system is a particular case of (2), corresponding to $M = N = P = 0$, the Chen system is obtained for $M = -a$, $N = 1 + c$, $P = 0$, the Lu system corresponds to $M = -c$, $N = 1 + c$, $P = 0$ and for $M = -a$, $N = 1$, $P = 1 - a$ is obtained the T-system.

In this paper we study the bifurcation of the equilibrium points and the existence of homoclinic and heteroclinic orbits of the general system (2). In Section 2 it is shown that the feedback control can be used in order to generate chaos in the stable Lorenz system. In Section 3 the existence of the heteroclinic orbits is correlated with the regular behavior of (2) and necessary conditions for the anticontrol of chaos are presented. Some results contained in [1, 5, 7, 10] are obtained as particular cases of our calculation.

## 2. EQUILIBRIUM POINTS

The system (2) is symmetric with respect to the $Oz$ axis. The dynamics of the system is affected by this symmetry: if $E_+(x_+, y_+, z_+)$, $E_-(-x_+, -y_+, z_+) = S_{Oz}(E_+)$ are equilibrium points then they are on the same type and $W^s(E_-) = S_{Oz}(W^s(E_+))$, respectively $W^u(E_-) = S_{Oz}(W^u(E_+))$; the stable (respectively unstable) manifolds of an equilibrium point situated on the $Oz$ axis are $S_{Oz}$ invariant (symmetric with respect to the $Oz$ axis). Besides, this hints that pithfork bifurcation for equilibria and periodic solutions is possible.

**Theorem 3.1.** *The system (LC) has the equilibrium point $O(0, 0, 0)$ for all values of the parameters a, b, c, M, N, P. It is the unique equilibrium situated in the Oz axis.*

*1. $O(0, 0, 0)$ is a saddle with $\dim(W^s) = 2$ and $\dim(W^u) = 1$ iff $b > 0$ and $a(M + N + c - 1) > 0$.*

*2. $O(0, 0, 0)$ is a saddle with $\dim(W^s) = 1$ and $\dim(W^u) = 2$ iff $b > 0$, $a(M + N + c - 1) < 0$ and $N - a - 1 > 0$.*

*3. $O(0, 0, 0)$ is attractor iff $b > 0$, $a(M + N + c - 1) < 0$ and $N - a - 1 < 0$.*

*4. If $a(M + N + c - 1)(N - a - 1) = 0$ then $O(0, 0, 0)$ is not hyperbolic.*

*5. If $\dfrac{b(M + N + c - 1)}{1 - P} > 0$ the system has three isolated equilibrium points: $O(0, 0, 0)$, $E_\pm\left(\pm\sqrt{\dfrac{b(M + N + c - 1)}{1 - P}}, \pm\sqrt{\dfrac{b(M + N + c - 1)}{1 - P}}, \dfrac{(M + N + c - 1)}{1 - P}\right)$.*

The proof comes directly from the analysis of the eingenvalues of the Jacobian matrix of (2) in $O(0, 0, 0)$.

The existence of three unstable equilibria can be associated with the chaotic behavior of the system. The pithfork bifurcation that occurs when $M + N + c - 1 = 0$ opens the way for a cascade of tangent and period doubling bifurcations from which the chaos emerges.

In the conditions $a > 0$, $b > 0$, which are fully accepted, the bifurcation line for the Lorenz system is $c = 1$. If $0 < c < 1$ the Lorenz system is stable because the origin is a global attractor. In the system (2) the pichfork bifurcation occours for small values of the control parameters $M, N, P$, when $0 < M + N = 1 - c < 1$, hence the feedback controller can generate chaos in the stable Lorenz system.

For example, in the Chen system the pitchfork bifurcation is obtained when $0 < c = a/2 < 1$, respectively $0 < c = a < 1$ in the Lu and the T systems. These particular results can be found in [1, 5, 8, 10]. An extensive analysis of bifurcations in Chen system is presented in [3].

### 3. HETEROCLINIC AND HOMOCLINIC ORBITS

Homoclinic and heteroclinic orbits are important concepts in the study of the bifurcation of vector fields and chaos. Many chaotic characteristics of a complex system are related to the existence of homoclinic and heteroclinic orbits. In this section we study the existence of these kind of orbits for the generalized Lorenz system (2). In order to prove the existence of the heteroclinic orbits a Lyapunov like function can be used. Its main properties are presented in the following lemma.

**Lemma 3.1.** *We consider the system (2) with $a > 0$, $b > 0$, $\frac{b-2a}{1-P} \geq 0$ and $N - 1 - a \leq 0$. We define*

$$V(x,y,z) = \frac{b(b-2a)}{1-P}(x-y)^2 + \left(z - \frac{x^2}{b}\right)^2 + \frac{b-2a}{2a}\left(x^2 - b\frac{(M+N+c-1)}{1-P}\right)^2.$$

*Then*

*a)* $\frac{dV}{dt}(\varphi(t, u_0)) \leq 0$ *for all* $t \in R$ *and* $u_0 \in R^3$, *hence* $\lim_{t \to \infty} V(\varphi(t, u_0)) \leq V(\varphi(t_2, u_0)) \leq V(\varphi(t_1, u_0)) \leq \lim_{t \to -\infty} V(\varphi(t, u_0))$ *for all* $t_1, t_2 \in R$, $t_1 < t_2$.

*b)* *if there exist* $t_1 < t_2$ *such that* $V(\phi(t_1, u_0)) = V(\phi(t_2, u_0))$ *then* $u_0$ *is an equilibrium point of the system.*

*c)* *If* $\lim_{t \to -\infty} \phi(t, u_0) = (0, 0, 0)$ *and there is* $t_3 \in R$ *such that* $x(0, u_0) > 0$ *then* $V(0, 0, 0) > V(\phi(t, u_0))$ *and* $x(t, u_0) > 0$ *for all* $t \in R$.

**Proof.** a) $\frac{dV}{dt}(\varphi)(t, u_0) = -2\frac{b(b-2a)(a+1-N)}{1-P}(x-y)^2 - 2b(bz - x^2)^2$

hence, in the hypothesis' conditions we obtain $\frac{dV}{dt}(\varphi(t, u_0)) \leq 0$. The inequality $\lim_{t \to \infty} V(\varphi(t, u_0)) \leq V(\varphi(t_2, u_0)) \leq V(\varphi(t_1, u_0)) \leq \lim_{t \to -\infty} V(\varphi(t, u_0))$ for all $t_1, t_2 \in R$, $t_1 < t_2$ is obvious because $V(\varphi(., u_0))$ is a decreasing function.

b) If $V(\phi(t_1, u_0)) = V(\phi(t_2, u_0)) = \alpha$ it results that $V(\phi(t, u_0)) = \alpha$ for all $t \in [t_1, t_2]$. In this case $\frac{dV}{dt}(\varphi(t, u_0)) = 0$ for all $t \in [t_1, t_2]$, hence $x(t) - y(t) = 0$ and $bz(t) - x^2(t) = 0$ for all $t \in [t_1, t_2]$. From (6) one can find that $x'(t) = y'(t) = z'(t)$ for all $t \in [t_1, t_2]$ hence $u_0$ is an equilibrium point.

c) The condition $\lim_{t \to -\infty} \phi(t, u_0) = (0, 0, 0)$ shows that $V(\phi(t, u_0)) \leq V(0, 0, 0)$ for all $t \in R$. Suppose that there is $t_2 \in R$ such that $V(\phi(t_2, u_0)) = V(0, 0, 0)$. For any $t_1 < t_2$ we have $V(0, 0, 0) = V(\phi(t_2, u_0)) \leq V(\phi(t_1, u_0)) \leq V(0, 0, 0)$. From b) it results that $V(\phi(t_2, u_0)) = V(\phi(t_1, u_0))$ and from a) it results that $u_0$ is an equilibrium point. Due to the fact that $x(t_3, u_0) > 0$ it results that $u_0 \neq 0$. This contradicts the hypothesis $\lim_{t \to -\infty} \phi(t, u_0) = (0, 0, 0)$, hence the supposition is not true.

Suppose that there is $t \in R$ such that $x(t, u_0) = 0$. Because

$$\frac{b(b-2a)}{1-P} y^2(t, u_0) + z^2(t, u_0) + \frac{b-2a}{2a}\left(b\frac{(M+N+c-1)}{1-P}\right)^2 \geq$$

$$\geq \frac{b-2a}{2a}\left(b\frac{(M+N+c-1)}{1-P}\right)^2$$

for all $t \in R$ it results that $V(0, y(t, u_0), z(t, u_0)) \geq V(0, 0, 0)$, which means that $u_0 = (0, 0, 0)$, contradiction with the hypothesis.

**Theorem 3.1.** *The system (2) with $a > 0$, $b > 0$, $\frac{b-2a}{1-P} \geq 0$ and $N - 1 - a \leq 0$ has the following properties*
 *a) the system has not closed orbits.*
 *b) the system (2) has not homoclinic orbits.*
 *c) if $b - 2a \geq 0$, $P < 1$ then $\lim_{t \to \infty} \phi(t, u_0) \in E(LC)$ for all $u_0 \in R^3$.*
 *d) if $b - 2a \geq 0$, $P < 1$, $c + M > 0$ and $\frac{(M+N+c-1)}{1-P} > 0$ then the system has two heteroclinic orbits that connect $(0, 0, 0)$ with the other equilibrium points. These orbits are two by two symmetric with respect to the Oz axis.*

**Proof.** a) Suppose that there is a closed orbit $O(u_0)$. It means that there are $t_1 < t_2$ such that $\phi(t_1, u_0) = \phi(t_2, u_0)$. In this conditions, for the function $V$ considered in Lemma 3.1., we have $V(\phi(t_1, u_0)) = V(\phi(t_2, u_0))$ which means that $u_0$ is an equilibrium point.

b) An orbit starting from $u$ is homoclinic to the equilibrium point $u_0$ if $u_0 = \lim_{t\to\infty} \varphi(t,u) = \lim_{t\to-\infty} \varphi(t,u)$. For the function $V$ considered in Lemma 3.1. we have $V(u_0) \leq V(\varphi(t,u)) \leq V(u_0)$ for all $t \in R$ which means that $u$ is an equilibrium point.

c) For the function $V$ considered in the previous theorem we have $0 \leq V^*(u_0) = \lim_{t\to\infty} V(\varphi(t,u_0)) \leq V(\varphi(t_1,u_0)) \leq V(u_0) \leq \lim_{t\to-\infty} V(\varphi(t,u_0))$ for all $t_1 \geq 0$. From the definition of $V$ it results that the orbit $O_+(u_0) = \{\varphi(t,u_0) | t \geq 0\}$ is bounded in $R^3$. Denote by $\Omega(u_0)$ the $\omega$-limit set of $O_+(u_0)$. For every $u \in \Omega(u_0)$ there is a sequence $(t_n)_{n\in N}$ of positive numbers such that $\lim_{t\to\infty} \varphi(t_n,u_0) = u$. Because $V$ is continuous it results that $V(u) = \lim_{t\to\infty} V(\varphi(t_n,u_0)) = V^*(u_0)$. Let consider $u_1 \in \Omega(u_0)$. Then $\varphi(t,u_1) \in \Omega(u_0)$ for all $t \in R$, hence $V^*(u_0) = V(u_1) = V(\varphi(t,u_1))$ for all $t \in R$. From Lemma 3.1. it results that $u_1$ is an equilibrium point of (2). It means that $\Omega(u_0) \subset E(S)$. But the set $\Omega(u_0)$ must be connected and there are at most three equilibrium points of (2), hence $\Omega(u_0)$ is formed by a single equilibrium point, *i.e.* $\Omega(u_0) = \lim_{t\to\infty} V(\varphi(t,u_0)) \in E(S)$.

d) In the conditions of d) the equilibrium point $O(0, 0, 0)$ is a saddle point with $\dim(W_O^s) = 2$ and $\dim(W_O^u) = 1$. There are also two equilibrium points $E_+$ and $E_-$ that are not situated on the $Oz$ axis. From Lemma 3.1. it results that the unstable manifold of $O$, denoted by $\Gamma_O$ is splitted in two parts that are $\varphi$ invariant:

$$\Gamma_O^+ = \left\{u_0 = (x_0, y_0, z_0), x_0 > 0, \lim_{t\to-\infty} \varphi(t,u_0) = O\right\}$$

$$\Gamma_O^- = \left\{u_0 = (x_0, y_0, z_0), x_0 < 0, \lim_{t\to-\infty} \varphi(t,u_0) = O\right\}$$

Then $\Gamma_O^+$ is a heteroclinic orbit joining $O$ with $E_+$ and $\Gamma_O^-$ is a heteroclinic orbit joining $O$ with $E_-$ and there are not other heteroclinic orbits. Because $\Gamma_O^- = S_{Oz}(\Gamma_O^+)$ it results that the heteroclinic orbits are symmetric with respect to the $Oz$ axis.

**Corollary 3.1.** *If $c > 1$ and $b \geq 2a > 0$ the Lorenz system has no closed or homoclinic trajectories. The system has three equilibrium and two heteroclinic orbits that connect $(0, 0, 0)$ with the other equilibrium points.*

**Corollary 3.2.** *If* $2c - a > 0$ *and* $(b - 2a)(c - a) \leq 0$ *the Chen system has no closed or homoclinic trajectories. The system has three equilibrium and two heteroclinic orbits that connect* (0, 0, 0) *with the other equilibrium points.*

Similar results, requesting stronger conditions, are presented in [7].

**Corollary 3.3.** *If* $c - a > 0$ *and* $b - 2a \leq 0$ *the T-system has no closed or homoclinic trajectories. The system has three equilibrium and two heteroclinic orbits that connect* (0, 0, 0) *with the other equilibrium points*.

This result, a particular case of our calculation, is reported in [10].

Remark: In the conditions of the previous theorem the behavior of (2) is very simple, every orbit of the system being attracted by an equilibrium and the chaotic behavior is excluded. In order to have chaos in (2) we must impose the condition $b < 2a$.